\documentclass{amsart}
\usepackage[english]{babel}
\usepackage{latexsym,bm}
\begin{document}
\newtheorem*{axiom}{Axiom}
\newtheorem*{definition}{Definition}
\newtheorem*{theorem}{Theorem}
\pagestyle{plain}
\title{Very Basic Set theory}
\author{Doeko H. Homan}
\date{January 22, 2026}
\begin{abstract}
The {\em naive conception of set} is described as  "Any predicate has an
extension". Then Russell's paradox is used to prove the inconsistency of naive
set theory. If the only constituent of a set is that set itself we call that 
set an `individual'. We show Russell's paradox permits sets that are
individuals thus that is not an inconsistency. Zermelo's axiomatization also
permits individuals. {\em Very basic set theory} is a set theory with
individuals and has a philosophical foundation in Ludwig Wittgenstein's
Tractatus Logico-Philosophicus.
\end{abstract}
\maketitle
\pagenumbering{arabic}
\section*{Introduction}\noindent
Georg Cantor defined a set as (Boolos [1971] p.215)
\begin{itemize}
\item "a `many, which can be thought of as one, i.e. a totality of definite\\
\mbox{ } elements that can be combined into a whole by a law'."
\end{itemize}
That means a set is `determined' by its elements. Sets are denoted by $s$, $t$,
$u$, \ldots . Symbol $\in$ is the `membership relation' symbol. Then $s\in t$
means $s$ is a `constituent' (also called `an element' or `a member') of the
totality $t$. The negation of $s\in t$ is denoted by $s\notin t$. Set $s$ is
an `individual' if the only constituent of $s$ is $s$ itself. A set with only
$s$ and $t$ as constituents is denoted by $\{s,t\}$ (Zermelo [1908] p.263).
Often the terminology is: `set $\{s,t\}$ contains only $s$ and $t$'.  But the
totality $\{s,t\}$ is `thought of as one' thus a set is a mental picture in
your mind. Therefore we prefer to write `constituent of a set', not `is
contained in a set'.\\
\mbox{ }\\
In (Oliver and Smiley [2006] $\S$5.1) Alex Oliver and Timothy Smiley explain
why\par
"numerous textbooks" prohibit individuals: individuals are not needed for\par
\mbox{ }mathematical purposes.\\
In (Kanamori [2003] pp.285-288) Akihiro Kanamori describes that with the
`empty\par
\mbox{ }extension' and the `axiom of restriction' there is no need for
individuals (also\par
\mbox{ }known as atoms or urelemente). Individuals "are objects distinct from
the\par
\mbox{ }null set yet having no members and capable of belonging to sets".\\
In (Boolos [1971] pp.216-217) George Boolos notes "individuals (= non-sets)".\\
In section 4.6.5 (Mendelson [2015]) Elliott Mendelson explains\par
"The words `individual' and `atom' are sometimes used as synonyms for\par
\mbox{ }`urelement'. Thus, urelements have no members". And in a footnote
(p.304)\par
"Zermelo's 1908 axiomatization permitted urelements".\\
In (Zermelo [1908] pp.262-265) Ernst Zermelo describes in "10. Theorem."\par
"F\"ur jedes Element $x$ von $M$ ist es definit, ob $x$ $\varepsilon$ $x$ oder
nicht; diese M\"oglichkeit\par
\mbox{ }$x$ $\varepsilon$ $x$ ist an und f\"ur sich durch unsere Axiome nicht
ausgeschlossen."\\
Thus `$x$ $\varepsilon$ $x$ or not' is definite; the possibility "$x$
$\varepsilon$ $x$" is not excluded by our axioms.\\
Therefore Zermelo's 1908 axiomatization permits sets `containing'
themselves.\par
The membership relation is definite, that means
\begin{itemize}
\item for every $s$ and for every $t$ applies either $s\in t$ or $s\notin t$,
but not both.
\end{itemize}\noindent
Then for every set $s$ it is the case $s\in s$ or $s\notin s$, but not both.
Therefore $s\in s$ is possible. But a set that `contains' itself is considered
problematic (Boolos [1971] p.217, pp.219-220)\par
"perhaps the mind ought to boggle at the idea of something's {\em containing}
itself."\par
"It is important to realize how odd the idea of something's containing itself
is."\par
"Here are some things. Now we bind them up into a whole. {\em Now} we have a
set."\\
However, `binding up some things' or `containing something' is not the
definition of a set. A set does not `contain' something, the constituents of a
set can be thought of as one totality.\par
Russell's paradox reads (Boolos [1971] p.217)\par
"No set can contain all and only those sets which do not contain themselves."\\
The negation of Russell's paradox reads
\begin{itemize}
\item for every $s$ exists a set $u$ such that (($u\in s$ and $u\in u$) or ($u
\notin s$ and $u\notin u$)).
\end{itemize}
If we choose for set $u$ set $s$ then the formula reads `for every $s$ ($s\in
s$ or $s\notin s$)'. It is the case that `for every $s$ ($s\in s$)' is
definite. Therefore Russell's paradox permits sets `containing' themselves.\\
\mbox{ }\par
In $\S$1 we describe sets and individuals. Atomic facts and elementary
propositions are defined which gives the foundation of set theory in the {\em
Tractatus} (Wittgenstein [1999]). In $\S$2 we prove the formulas of set theory
preserve definiteness. In $\S$3 we define equality for sets. Then axioms
specification, pairs and union are formulated. We prove regularity and some
other important theorems. In $\S$7 natural numbers are defined.
\section{What are individuals and sets?}
At an early age you develop the mental picture of `what is a set'. You belong
to a family or a clan, to the inhabitants of a village or a particular region.
Then experience shows there are things that constitute a set. The members of a
set are the constituents of that set.\par
Set theory is a language that follows rules of logic to communicate the mental
picture of a set with the physical world. Then someone reading, seeing or
listening that language forms a mental picture of a set. Set theory is made up
of `formulas'. A formula is a proposition or predicate about the facts of sets.
For every set $s$ and for every set $t$ the fact is the definiteness of $s\in
t$. A formula which is not complex or not a combination of other formulas is an
`atomic fact'. Formula $s\in t$ is an atomic fact, formula $t\in s$ is an
atomic fact. Formula $s\notin t$ is the negation of $s\in t$ thus $s\notin t$
is a complex formula therefore $s\notin t$ is not an atomic fact. Atomic facts
are what makes formulas true or false.
\begin{definition}
Set $s$ is an individual if the only constituent of set $s$ is set $s$ itself.
\end{definition}
A formula presupposes propositions or predicates while propositions or
predicates presupposes atomic facts. Then if an atomic fact is analyzed as
fully as possible the constituents finally reached are individuals. Individuals
cannot be analyzed further. In this way the naming of individuals is what is
logically first in logic (Wittgenstein [1999] prop. 3.26) and Bertrand Russell
(Introduction p.12).\par
If $s$ is a set then $s$ is not a proposition or predicate. If $s$ is an
individual then $s\in s$ therefore $s$ is the 'propositional connexion'
(Wittgenstein [1999] prop. 4.221), (Anscombe [1963] p.29-31). Thus individuals
are `elementary propositions'. Therefore in {\em Very basic set theory} the
naming of individuals is logically first.\par
\section{Logic and definiteness}
In this section we prove logical formulas preserve definiteness. Formulas are
denoted by $\Phi$, $\Psi$, $\ldots$ and are built up from formulas by means of
negation $\neg\Phi$ (not $\Phi$) and conjunction $\Phi\wedge\Psi$ ($\Phi$ and
$\Psi$). Negation takes precedence over conjunction. If $\Phi$ is a formula
then $(\Phi)$ is a formula. The parenthesis $($ and $)$ preclude ambiguity and
improve legibility. Then $s\notin t$ is the formula $\neg(s\in t)$.\par
If $\Phi$ is a formula then $\forall s(\Phi)$ and $\exists s(\Phi)$ are
formulas. $\forall s(\Phi)$ means `for every set it is the case $\Phi$'. The
negation of $\forall s(\neg\Phi)$ is denoted by $\exists s(\Phi)$ thus $\exists
s(\Phi)$ means `exists a set such that $\Phi$'. The parenthesis following
$\forall s$ and $\exists s$ are mandatory. However, e.g. without ambiguity
$\forall s\exists t(s\in t)$ is short for the formula $\forall s(\exists t(s\in
t))$.\par
If there is $s$ in $\Phi$ not guided in $\Phi$ by $\forall s$ or by $\exists s$
then `$s$ is free in $\Phi$'. Thus $s$ is free and $t$ is not free in formula
$\forall t(s\in t)$.\par
Set $s$ in $\forall s$ or in $\exists s$ is not a `variable'. $\forall s$ means
`for every set', thus also for a set denoted by $s$. And $\exists s$ means
`exists a set' denoted by $s$. Set theory is a {\em language}.\par
We define the logical connectives $\vee$, $\rightarrow$ and $\leftrightarrow$.
\begin{description}
\item{} $\Phi\vee\Psi$ ($\Phi$ or $\Psi$ or both) is short for the formula
$\neg(\neg\Phi\wedge\neg\Psi)$,
\item{} $\Phi\rightarrow\Psi$ ($\Phi$ implies $\Psi$) is short for the formula
$(\neg\Phi\vee\Psi)$,
\item{} $\Phi\leftrightarrow\Psi$ ($\Phi$ if and only if $\Psi$) is short for
the formula ($(\Phi\rightarrow\Psi)\wedge(\Psi\rightarrow\Phi))$.
\end{description}
There are many rewrite rules: $\neg(\neg\Phi)$ is logically equivalent to
$\Phi$, $\neg(\Phi\leftrightarrow\Psi)$ is logically equivalent to $(\Phi
\leftrightarrow\neg\Psi)$ and $((\Phi\wedge\Psi)\wedge\neg(\Phi\wedge\Psi))$
is logically equivalent to $((\Phi\wedge\Psi)\wedge(\neg\Phi\vee\neg\Psi))$.
More rewrite rules in (Mendelson [2015] p.15).
\begin{definition} Formula $\Phi$ is `definite' if it is not the case that both
$\Phi$ and $\neg\Phi$.
\hfill\boldmath$\Box$\end{definition}
\begin{theorem} The building up of formulas preserve definiteness.
\end{theorem}\noindent
{\em Proof.} For every $s$ and for every $t$ it is not the case that both $s
\in t$ and $\neg(s\in t)$\par
\mbox{ } therefore atomic facts and elementary propositions are definite.\par
If $\Phi$ is definite then it is not the case that both $\Phi$ and $\neg\Phi$.
Thus it is not the\par
\mbox{ } case that both $\neg(\neg\Phi)$ and $\neg\Phi$ therefore $\neg\Phi$ is
definite.\par
If $\Phi$ is definite then $\forall s(\Phi)$ and $\exists s(\Phi)$ are definite
otherwise $(\forall s(\Phi)\wedge\neg\forall s(\Phi))$ or\par
\mbox{ } $(\exists s(\Phi)\wedge\neg\exists s(\Phi))$ thus $\exists s(\Phi
\wedge\neg\Phi)$.\par
And $(\Phi\wedge\Psi)\wedge\neg(\Phi\wedge\Psi)$ is logically equivalent to\par
\mbox{ } $((\Phi\wedge\neg\Phi)\wedge\Psi)\vee(\Phi\wedge(\Psi\wedge\neg\Psi))$
therefore if $\Phi$ and $\Psi$ are definite then\par
\mbox{ } $(\Phi\wedge\Psi)$ is definite.\par
The logical connectives `or', `implies' and `if and only if' are expressed in
terms\par
\mbox{ } of logical `not' and logical `and'. Thus if $\Phi$ and $\Psi$ are
definite then $(\Phi\vee\Psi)$,\par
\mbox{ } $(\Phi\rightarrow\Psi)$ and $(\Phi\leftrightarrow\Psi)$ are
definite.\par
If an atomic fact is analyzed as fully as possible constituents finally reached
are\par
individuals. Atomic facts and elementary propositions are definite.\par
Therefore the building up of formulas preserve definiteness.
\hfill{\boldmath$\Box$}
\section{Equality for sets}
Experience shows a set `is equal to' or `is not equal to' another set. A set is
determined by its constituents. Thus sets with different constituents are
different sets, and the constituents of different sets are different. We define
equality for sets.
\begin{definition} Set $s$ `is equal to' set $t$, denoted by $s=t$ if $\forall
u(u\in s\leftrightarrow u\in t)$.\par
Set $s$ `is not equal to' set $t$, denoted by $s\neq t$ if $\neg(s=t)$.
\hfill\boldmath$\Box$\end{definition}\noindent
The axioms for equality are
\begin{itemize}
\item \mbox{ }\mbox{ }(equality) $\forall s\forall t(s=t\rightarrow\forall u(s
\in u\leftrightarrow t\in u))$,
\item (individuals) $\forall s(s\in s\rightarrow\forall u(u\neq s\rightarrow
u\notin s))$.
\end{itemize}
Thus if $s=t$ then $s$ and $t$ have the same constituents and $s$ and $t$ are
members of the same sets. And the only constituent of an individual is that
individual itself.\\
Then $\forall s\forall t((s\in s\wedge t\notin t)\rightarrow(s\neq t\wedge t
\notin s))$.
\section{Axiom `specification'}
If $\Phi$ is definite and if exists set $u\in s$ such that $\Phi$ then it is
not the case that for every set $u$ $(u\in s\rightarrow\neg\Phi)$. Thus we can
postulate a set $v$ whose constituents are all and only those constituents of
$s$ such that for every constituent of $v$ it is the case $\Phi$. That is axiom
specification (also known as `separation').
\begin{itemize}
\item if $u$ is free in $\Phi$ and $\Phi$ is definite then\\
$\forall s(\exists u(u\in s\wedge\Phi)\rightarrow\exists v\forall u(u\in v
\leftrightarrow(u\in s\wedge\Phi)))$.
\end{itemize}
Thus the conjunction of mutually exclusive properties does not constitute a
set.
\begin{theorem}
If an individual is a constituent of set $s$ then exists set $v$ whose\par
constituents are all and only those individuals that are constituent of $s$.
\end{theorem}\noindent
{\em Proof.} Formula $(u\in s\wedge u\in u)$ is definite and is not the
conjunction of mutually\par
exclusive properties. Apply specification to set $s$ to find set $v$\par
\mbox{ }$\exists v\forall u(u\in v\leftrightarrow(u\in s\wedge u\in u))$. Then
an individual is a constituent of $v$.\par
The formula applies to every $u$ thus we can choose $v$\par
\mbox{ }$\exists v(v\in v\leftrightarrow(v\in s\wedge v\in v))$ therefore\par
\mbox{ }$\exists v((v\notin v\vee (v\in s\wedge v\in v))\wedge ((v\notin s\vee
v\notin v)\vee v\in v))$.\par
Then $\exists v(v\notin v\vee v\in s)$.\par
Therefore $v$ is not an individual, or $v$ is an individual and a constituent
of $s$.
\hfill{\boldmath$\Box$}
\begin{theorem} The set of all sets is nonexisting.
\end{theorem}\noindent
{\em Proof.} If $\exists u(u\in s\wedge u\notin u)$ then formula $(u\in s\wedge
u\notin u)$ is definite and is not the\par
conjuction of mutually exclusive properties.\par
Apply specification to set $s$ to find set $v$\par
\mbox{ }$\exists v\forall u(u\in v\leftrightarrow(u\in s\wedge u\notin u))$.
Then a set is a constituent of $v$.\par
The formula applies to every $u$ thus we can choose $v$\par
\mbox{ }$\exists v(v\in v\leftrightarrow(v\in s\wedge v\notin v))$
therefore\par
\mbox{ }$\exists v((v\notin v\vee(v\in s\wedge v\notin v))\wedge((v\notin s\vee
v\in v)\vee v\in v))$.\par
Then $\exists v(v\notin v\wedge v\notin s)$.\par
Thus $v$ is not an individual and $v$ is not a constituent of $s$.\par
Therefore the set of all sets is nonexisting.
\hfill{\boldmath$\Box$}
\section{Axiom `pairs'}
Any pair of sets constitute a set. We postulate `pairs'.
\begin{itemize}
\item (pairs) $\forall s\forall t\exists v\forall u(u\in v\leftrightarrow(u=s
\vee u=t))$.
\end{itemize}
Set $v$ is denoted by $\{s,t\}$, $s$ and $t$ are the only constituents of
$\{s,t\}$.\\
Then $\{s,t\}=\{t,s\}$ thus $\{s,t\}$ is also called `an unordered pair'.\\
The `singleton' $\{s\}$ of set $s$ is $\{s,s\}$. Then\par
\mbox{ } $s\notin s\rightarrow\{s\}\neq s$ otherwise $s\in s$. Therefore $s
\notin s\rightarrow\{s\}\notin\{s\}$,\par
\mbox{ } $s\in s\rightarrow\{s\}=s$, the singleton of an individual is that
individual itself,\par
\mbox{ } $s\neq t\rightarrow\{s,t\}\notin\{s,t\}$. Thus a pair of different
sets is not an individual.\\
However, if $s\neq t$ and $(s\in t\wedge t\in s)$ then $s\notin s$ and $t\notin
t$. Then exists set $\{s,t\}$ and the `vicious membership circle' $(s\in t
\wedge t\in s)$. That is contradictory to `the naming of individuals is
logically first'. Therefore if $s\neq t$ it is not the case that\par
$\forall v(v\in\{s,t\}\rightarrow\exists u((u\in\{s,t\}\wedge u\in v)\wedge
u\notin u))$.\\
The naming of individuals is logically first. Thus vicious membership circles
do not exist. We prove `regularity'.
\begin{theorem} It is not the case that exists set $w$ such that\par
$\forall v(v\in w\rightarrow\exists u((u\in w\wedge u\in v)\wedge u\notin
u))$.
\end{theorem}\noindent
{\em Proof.}
A finite sequence of sets \mbox{ }$s$, $t$, $u$, \ldots, $z$ \mbox{ }such that
\mbox{ }$s\in t$, $t\in u$, \ldots, $z\in s$\par
is nonexisting otherwise there are vicious membership circles. And an
`infinite'\par
sequence of sets \mbox{ }$s$, $t$, $u$, $v$, \ldots \mbox{ }such that each term
is a member of the previous\par
one \mbox{ }$t\in s$, $u\in t$, $v\in u$ \ldots \mbox{ }is nonexisting
otherwise for every set $x$ belonging to\par
the sequence exists set $y$ belonging to that sequence such that $y\in x$.\par
That is contradictory to `the naming of individuals is logically first'.
\hfill{\boldmath$\Box$}\\
\mbox{ }\\
Formula $\neg\exists w\forall v(v\in w\rightarrow\exists u((u\in w\wedge u\in
v)\wedge u\notin u))$ is definite. Therefore
\begin{itemize}
\item (regularity) $\forall w\exists v(v\in w\wedge\forall u((u\in w\wedge u\in
v)\rightarrow u\in u))$.
\end{itemize}
The formula applies to every $u$ thus we can choose $v$. Then the formula
reads\par
$\forall w\exists v(v\in w\wedge((v\notin w\vee v\notin v)\vee v\in v))$, that
is $\forall w\exists v(v\in w)$.\\
Therefore an `empty set' is nonexisting. If an individual is a constituent of
$w$ then $v$ is an individual. And if $\exists v(v\notin v\wedge v\in w)$ then
$\exists v(v\notin v\wedge v\in w)$ such that $v$ and $w$ only have individuals
in common.
\section{Axiom `union'}
The constituents of the members of a set constitute a set. We postulate
`union'.
\begin{itemize}
\item (union) $\forall s\exists v\forall u(u\in v\leftrightarrow\exists t(t\in
s\wedge u\in t))$.
\end{itemize}
Set $v$ is denoted by $\bigcup s$. If $s$ is an individual or a set of
individuals then $\bigcup s=s$. The `union $s\cup t$' is defined by $s\cup t=
\bigcup\{s,t\}$. Then $\bigcup\{\{s\},\{t\}\}=\{s,t\}$ and $s\cup(t\cup u)=
(s\cup t)\cup u$. If $s$ and $t$ are individuals then $s\cup t=\{s,t\}$.
\begin{theorem} Any sequence of different sets \mbox{ }$s$, $\bigcup s$,
$\bigcup(\bigcup s)$, \ldots \mbox{ }is finite.
\end{theorem}\noindent
{\em Proof.} Then $s\neq\bigcup s\rightarrow\exists t(t\in s\wedge t\notin t)$
otherwise $s=\bigcup s$.\par
$\bigcup s\neq\bigcup(\bigcup s)\rightarrow\exists u\exists t(u\in t\wedge t\in
s\wedge u\notin u\wedge t\notin t)$ otherwise $\bigcup s=\bigcup
(\bigcup s)$,\par
$\bigcup(\bigcup s)\neq\bigcup(\bigcup(\bigcup s))\rightarrow\exists v\exists u
\exists t(v\in u\wedge u\in t\wedge t\in s\wedge v\notin v\wedge u\notin u
\wedge t\notin t)$,\par
\mbox{ }\mbox{ }\mbox{ } \ldots \mbox{ }.\par
Thus there is a rule to construct for every set $z\neq s$ belonging to the
sequence a\par
series of sets $t$, $u$, $v$, \mbox{ }\ldots\mbox{ }, $x$, $y$, $z$ such that
each term is a member of the previous\par
one: $t\in s$, $u\in t$, $v\in u$, \mbox{ }\ldots\mbox{ }, $x\in v$, $y\in x$,
$z\in y$ contradictory to regularity if\par
the sequence is `infinite': the naming of individuals is logically first.\par
Therefore the sequence terminates. Thus exists set $w$ belonging to the
sequence\par
such that $\bigcup w=w$. Then $s\cup\bigcup s\cup\bigcup(\bigcup s)\cup\ldots$
$w$ is a set, and is called `the\par
transitive closure of $s$'.
\hfill{\boldmath$\Box$}\\
\mbox{ }\\
The definitions of equality and axiom pairs are only meaningful if exist sets
the naming of which are logically first. Thus with individuals and regularity.
And regularity is used with axiom union to prove that for every set the
transitive closure exists. The proof requires no other axioms or no `set of all
natural numbers'.
\section{Natural numbers}
In (Halmos [1974] p.47) the Peano Axioms for the set of all natural numbers
"used to be considered as the fountainhead of all mathematical knowledge". In
this section we define `transitive sets' and `number sequences' satisfying the
Peano Axioms for natural numbers.
\begin{definition}
Set $s$ is a `subset' of $t$, denoted by $s\subseteq t$, if $\forall u(u\in s
\rightarrow u\in t)$.\par
$s$ is transitive if $\forall u(u\in s\rightarrow u\subseteq s)$, or
equivalently $\bigcup s\subseteq s$.
\hfill\boldmath$\Box$\end{definition}\noindent
Every individual and every set of individuals is transitive. If $s$ is
transitive then $\bigcup s$ and $s\cup\{s\}$ are transitive.
\begin{theorem} If $s$ is transitive and $\exists u(u\in s\wedge u\notin u)$
then a set of individuals is a\par
member of $s$: $\exists v(v\in s\wedge v\notin v\wedge\forall u(u\in v
\rightarrow u\in u))$.
\end{theorem}\noindent
{\it Proof.} Apply regularity to $s$ to find $v\in s\wedge v\notin v$ such
that\par
$\forall u((u\in s\wedge u\in v)\rightarrow u\in u)$. Set $s$ is transitive and
$v\in s$ thus $\forall u(u\in v\rightarrow u\in s)$.\par
Therefore $\exists v(v\in s\wedge v\notin v\wedge\forall u(u\in v\rightarrow u
\in u))$.
\hfill{\boldmath$\Box$}\\
\mbox{ }\\
If the only individuals that are constituent of transitive set $s$ are $u$ and
$v$ then $\{u,v\}$ is the only possible set of individuals that is a member of
$s$ therefore $\{u,v\}\in s\cup\{s\}$.\par
If $s$ and $t$ are transitive sets with transitive members and the only
individuals that are members of $s$ or of $t$ are $u$ and $v$, thus $\forall w
((w\in s\cup t\wedge w\in w)\rightarrow w\in\{u,v\})$, then the `law of
trichotomy' for $s$ and $t$ is easily proved and reads
$$(s\notin s\wedge t\notin t)\rightarrow(s\in t\vee s=t\vee t\in s).$$
If $s$ is a transitive set with transitive members then $\bigcup s$ is a
transitive set with transitive members and $(s\notin s\rightarrow\bigcup s
\notin\bigcup s)$. And $s\notin\bigcup s$ otherwise $s\in s$. Therefore
$\bigcup s=s$ or $s=\bigcup s\cup\{\bigcup s\}$.\par
It is straightforward to define a number sequence satisfying the Peano Axioms
for natural numbers. The first axiom of the Peano Axioms reads `there is a
first natural number'. Natural numbers are also known as `finite ordinals'.
\begin{definition}
Set $s$ is a `natural number with first number $\alpha$' if $s$ is a
transitive\par
set with transitive members, $\alpha$ is a pair of individuals and
$$\bigcup\alpha=\alpha\wedge s\notin\alpha\wedge\forall u(u\in s\cup\{s\}
\rightarrow(u\in\alpha\cup\{\alpha\}\vee\bigcup u\in u)).$$\par
The `successor' of natural number $s$ is $s\cup\{s\}$.
\hfill\boldmath$\Box$\end{definition}\noindent
Then $\alpha$ is a first natural number. In physical reality a first number is
called `the unit of measurement'. The classical Greek concept of number (\em
arithmos\em) is "a finite plurality composed of units, where a unit is whatever
counts as one thing in the number under consideration." (Mayberry [2000]
p.21).\par
If $z$ (`zero'), $x$ and $y$ are individuals then there is number sequence $S$
with first number $\{z,x\}$ (`x coordinate') and there is number sequence $T$
with first number $\{z,y\}$ (`y coordinate'). All elements of both sequences
have $z$ in common.\\
If $u$ belongs to the $S$ sequence then $x\in u$ implies $u$ is an $x$
coordinate, if $v$ belongs to the $T$ sequence then $y\in v$ implies $v$ is an
$y$ coordinate. Then $\{u,v\}$ is an `ordered pair' with an intrinsic intuitive
meaning. The Kuratowski ordered pair does not have any intrinsic intuitive
meaning. It is "just a convenient way" to define ordered pairs (Mendelson
[2015] p.235, Kanamori [2003] \S 5).\par
Any pair of individuals is a first number of a natural number sequence.
Therefore `the set of all natural numbers' is nonexisting. However, it is not
prohibited to postulate `set $\omega$ of all natural numbers with first number
$\alpha$'. Then $\bigcup\omega=\omega$. Therefore we can define a number
sequence with first number $\omega$ analogous to natural numbers with first
number $\alpha$. Then there are the finite ordinals (natural numbers) with
first number $\alpha$, followed by `transfinite ordinal numbers with first
number $\omega$'. And if we define $\omega+\omega=\omega\cdot2$ then $\bigcup
(\omega\cdot2)=\omega\cdot2$ thus the next sequence of transfinite ordinals are
the ordinals with first number $\omega\cdot2$, and so on, that is\par
$\alpha$, $\alpha+1$, $\alpha+2$, \ldots , $\omega$, $\omega+1$, $\omega+2$,
\ldots , $\omega\cdot2$, $\omega\cdot2+1$, $\omega\cdot2+2$, \ldots ,
$\omega\cdot3$, \ldots .\\
Thus the logic in {\em Very basic set theory} is capable to deal with
transfinite numbers.
\section*{Conclusion}
{\em Very basic set theory} is a set theory with the philosophical
foundation in Ludwig Wittgenstein's Tractatus Logico-Philosophicus. It is a set
theory with individuals. The atomic facts are the membership relations and the
elementary propositions are the individuals. The formulas of set theory
preserve definiteness. The naming of individuals is logically first and the
logic is capable to deal with transfinite numbers.\\
\mbox{ }\par
We list some differences with `mainstream set theory'.
\begin{itemize}
\item There is no `empty set' (Kanamori [2003] \S1).
\item The singleton of an individual is that individual itself (Kanamori [2003]
\S2).
\item Regularity is not an axiom but a theorem.
\item The proof that `every set has a transitive closure' requires no other
axioms or no `set of all natural numbers' (Mendelson [2015] pp.286-287).
\item Cantor's Theorem (Zermelo [1908] p.276) reads\\
\mbox{ }\mbox{ }for every set $s$, the `number of members' of $s$ is less than
the `number of\\
\mbox{ }\mbox{ }members' of the powerset ${\mathcal{P}}s$.\\
However, there is no `empty set' therefore ${\mathcal{P}}\{s\}=\{\{s\}\}$.\\
Thus Cantor's Theorem does not hold.
\item Any pair of individuals is the first number of a sequence of natural
numbers. Then transfinite ordinals are postulated. If e.g. $\bigcup\omega=
\omega$ then $\omega$ is the first number of a sequence of transfinite
ordinals. Therefore $\omega$ is not a `limit ordinal', $\omega$ is the first
number of a sequence satifying the Peano Axioms for natural numbers.
\end{itemize}
\newpage
\section*{References}
\begin{description}
\item G.E.M. Anscombe. {\em An Introduction to Wittgenstein's Tractatus}\\
Second Edition, Revised in 1963 (https://archive.org)\\
First HARPER TORCHBOOK edition published 1965 by\par
Harper \& Row, Publishers, Incorporated
\item George Boolos. {\em The iterative conception of set}\\
The Journal of Philosophy\\
Volume LXVIII, no. 8, April 22, 1971.
\item P.R. Halmos. {\em Naive Set Theory}\\
Springer-Verlag New York Inc., 1974.
\item Akihiro Kanamori. {\em The empty set, the singleton and the ordered
pair}\\
The Bulletin of Symbolic Logic, Volume 9, Number 3, Sept. 2003.
\item Elliott Mendelson. {\em Introduction to Mathematical Logic}\\
Sixth Edition, 2015.
\item J.P. Mayberry. {\em The Foundations of Mathematics in the Theory of
Sets}\\
Cambridge University Press, 2000.
\item Alex Oliver and Timothy Smiley. {\em What Are Sets and What Are They
For?}\\
Philosophical Perspectives, 2006, Vol. 20, Metaphysics (2006), pp. 123-155\\
Stable URL: https://www.jstor.org/stable/4494502
\item Ludwig Wittgenstein. {\em Logisch-philosophische Abhandlung}\\
Translated by C.K. Ogden. {\em Tractatus Logico-Philosophicus}\\
With an Introduction by Bertrand Russell, 1922\\
Dover Publications, Inc. 1999.
\item Ernst Zermelo. {\em Untersuchungen \"uber die Grundlagen der Mengenlehre
I.}\\
Mathematische Annalen (1908) Vol. 65 (https://eudml.org/doc/158344)
\end{description}
\end{document}